\documentclass[twocolumn]{autart}    

\usepackage{amsmath,amsfonts,amssymb,mathtools}
\usepackage{bbm}
\usepackage{graphicx}
\usepackage{exscale}
\usepackage{stackrel}

\usepackage{clrscode}
\usepackage{array}
\usepackage{hhline}
\usepackage{color, colortbl}
\definecolor{myGray}{rgb}{0.95,0.95,0.95}

\newtheorem{exampe}{Example}[section]

\begin{document}

\begin{frontmatter}
\runtitle{MATLAB-based Unscented Kalman Filtering}  

\title{Continuous-discrete unscented Kalman filtering framework by MATLAB ODE solvers and square-root methods} 

\thanks[footnoteinfo]{\scriptsize This paper was not presented at any IFAC
meeting. Corresponding author M.~V.~Kulikova.\\
The authors acknowledge the financial support of the Portuguese FCT~--- \emph{Funda\c{c}\~ao para a Ci\^encia e a Tecnologia}, through the projects UIDB/04621/2020 and UIDP/04621/2020 of CEMAT/IST-ID, Center for Computational and Stochastic Mathematics, Instituto Superior T\'ecnico, University of Lisbon.}

\author[IST]{Maria V. Kulikova}\ead{maria.kulikova@ist.utl.pt},
\author[IST]{Gennady Yu. Kulikov}\ead{gennady.kulikov@tecnico.ulisboa.pt}

\address[IST]{CEMAT, Instituto Superior T\'{e}cnico, Universidade de Lisboa,
          Av. Rovisco Pais 1,  1049-001 Lisboa, Portugal}


\begin{abstract}
This paper addresses the problem of designing the {\it continuous-discrete} unscented Kalman filter (UKF) implementation methods. More precisely, the aim is to propose the MATLAB-based UKF algorithms for {\it accurate} and {\it robust} state estimation of stochastic dynamic systems. The accuracy of the {\it continuous-discrete} nonlinear filters heavily depends on how the implementation method manages the discretization error arisen at the filter prediction step. We suggest the elegant and accurate implementation framework for tracking the hidden states by utilizing the MATLAB built-in numerical integration schemes developed for solving ordinary differential equations (ODEs). The accuracy is boosted by the discretization error control involved in all MATLAB ODE solvers. This keeps the discretization error below the tolerance value provided by users, automatically. Meanwhile, the robustness of the UKF filtering methods is examined in terms of the stability to roundoff. In contrast to the pseudo-square-root UKF implementations established in engineering literature, which are based on the one-rank Cholesky updates, we derive the stable square-root methods by utilizing the $J$-orthogonal transformations for calculating the Cholesky square-root factors.
\end{abstract}

\end{frontmatter}

\section{Introduction}\label{sect1}

The optimal Bayesian filtering for nonlinear dynamic stochastic systems is a well-established area in engineering literature. A wide range of the filtering  methods has been developed under assumption of Gaussian systems and by applying a variety of numerical methods for computing the multidimensional Gaussian integrals required for finding the mean and covariance to construct the suboptimal solution~\cite{Ito2000}. Among well-known Gaussian filters designed, we may mention the Gauss-Hermit quadrature filters (GHQF) in~\cite{Haykin2007}, the third- and high-degree Cubature Kalman filters (CKF) in~\cite{2010:Haykin,2013:Automatica:Jia,2018:Haykin}, respectively. Nevertheless, the Gaussian property is rarely preserved in nonlinear Kalman filtering (KF) realm and is often violated while solving real-world stochastic applications. In contrast to the Gaussian filters, the Unscented Kalman filter (UKF) calculates true mean and covariance even in non-Gaussian stochastic models~\cite{2000:Julier}. Motivated by this and some other benefits of the UKF framework, we suggest a few novel MATLAB-based UKF implementation methods for {\it accurate} and {\it robust} state estimation of stochastic dynamic systems. The accuracy of any {\it continuous-discrete} filtering algorithm heavily depends on how the implementation method manages the discretization error arisen. The robustness of the filtering method with the KF-like structure is examined with respect to roundoff errors and the feasibility to derive numerically stable square-root counterparts. Table~\ref{Tab:survay} presents a systematic classification of {\it continuous-discrete UKF implementation methods} already established in engineering literature as well as open problems to be solved in future. It provides a basis for understanding various implementation strategies depending on discretization accuracy and numerical stability to roundoff.

\begin{table*}[ht!]
{\tiny
\renewcommand{\arraystretch}{1.3}
\caption{A systematization of the continuous-discrete UKF implementation methods available at present. } \label{Tab:survay}
\centering
\begin{tabular}{|l|c|c|l|cc|c|l|}
\hline
& \multicolumn{2}{c|}{\bf Type~I. Discrete-Discrete Approach} & & \multicolumn{4}{l|}{\bf Type~II. Continuous-Discrete Approach}  \\
\hhline{~|-|-|~|-|-|-|-}
\hhline{~|-|-|~|-|-|-|-}
{\bf Unscented }& {\bf Euler-Maruyama} & {\bf It\^{o}-Taylor} & {\bf Key properties} & \multicolumn{2}{r|}{\bf MDE} & {\bf SPDE} &  \\
 {\bf  Kalman}          & \cite[Chapter~9]{1999:Kloeden:book} & \cite[Chapter~5]{1999:Kloeden:book} & {\bf of the examined}   & \multicolumn{2}{r|}{\cite[Equation~(34)]{2007:Sarkka}} & \cite[Eq.~(35)]{2007:Sarkka} &  \\
 {\bf filtering}    & {\tiny (strong order 0.5)} & {\tiny (strong order 1.5)}       &  {\bf filtering methods}        & \multicolumn{2}{r|}{\bf $n^2+n$ eqs.}  & {\bf $2n^2+n$ eqs.} & {\bf \tiny \quad Benefits} \\
\hline
\hhline{~|-|-|~|-|-|-|-}
{\bf Non-square-root} &  studied in~\cite{2019:Leth},  & studied in~\cite{lyons2014series},  & $\to$\,unstable\,to\,roundoff  & $\stackrel[{\tiny \rm{\;\;ODE\;solutions:}}]{ \bullet\,\rm{ fixed-stepsize}}{}$ & $\stackrel[\text{\makebox[0pt]{\cite{takeno2012numerical}}}]{\text{\makebox[0pt]{\cite{2007:Sarkka}}}}{}$ & \cite{2007:Sarkka} & $\stackrel[{\tiny \rm{to\,Type~I}}]{ \rm{\to comparable}}{\mathrm{}}$ \\
\cline{5-8}
{\bf implementations} &  \cite{lyons2014series} etc. &  etc. & \,and\,to\,discretization\,$\leftarrow$  & $\stackrel[{\bullet\,\rm{MATLAB\,ones:}}]{ \bullet\,\rm{NIRK\;solvers:\,\;}}{}$ & $\stackrel[{\bf this\,paper}]{\bf open}{}$
 & $\stackrel[\text{\makebox[0pt]{\cite{KuKu20dIFAC}}}]{\text{\makebox[0pt]{\cite{2017:SP:Kulikov}}}}{}$
& {\tiny $\stackrel[\to \rm{LE\,control}]{\to\rm{GE\,control}}{}$ }\\
\hline
{\bf Square-root\,ones}  &  \multicolumn{2}{c|}{$\to$\,ensures\,symmetry\,\&\,positivity\,of\,$P_{k|k}$}  & {$\to$ improved stability}$\leftarrow$  & \multicolumn{2}{l|}{{\bf Theoretically derived} } &  & $\Downarrow$ This means \\
{\bf $\bullet$ Cholesky\,type}  &                & & { $\bullet$\,triangular\,SR\,factors}        &    \multicolumn{2}{l|}{{\bf SR solution:} \cite[Eq.~(64)]{2007:Sarkka} } & \cite[Eq.~(35)]{2007:Sarkka}  & Methods with \\
-- 1 rank updates  &    {\bf open}                        & \cite{KuKu21aEJCON} & {$\to$`pseudo-SR'\,variants}$\leftarrow$                 &
$\stackrel[{\bullet\,\rm{MATLAB\,ones:}}]{ \bullet\,\rm{NIRK\;solvers:\,\;}}{}$ & $\stackrel[{\bf this\,paper}]{\bf open}{}$ & 
$\stackrel[{\bf this\,paper}]{\text{\makebox[0pt]{\cite{2020:ANM:KulikovKulikova}}}}{}$ & LE control are, \\
\cline{2-7}
-- JQR-based alg. &    \cite{KuKu20bIFAC}                          & \cite{KuKu21aEJCON} & {$\to$true\,SR\,methods}$\leftarrow$                    &
$\stackrel[{\bullet\,\rm{MATLAB\,ones:}}]{ \bullet\,\rm{NIRK\;solvers:\,\;}}{}$ & $\stackrel[{\bf this\,paper}]{\bf open}{}$ & $\stackrel[\text{\makebox[0pt]{\cite{KuKu20dIFAC}}}]{\text{\makebox[0pt]{\cite{2020:ANM:KulikovKulikova}}}}{}$ &   in general, \\
\cline{2-7}
{\bf $\bullet$ SVD\,type}  &   {\bf open}                & \cite{KuKu21IEEE_TAC} & { $\bullet$\,full\;matrix\,SR\,factors}       &   $\stackrel[{ \bullet\,\rm{MATLAB\,ones:}}]{ \bullet\,\rm{NIRK\;solvers:\,\;}}{}$ & $\stackrel[{\bf open}]{\bf open}{}$
 & $\stackrel[{\bf open}]{\bf open}{}$ & faster, but   \\
\hhline{-|-|-|~|-|-|-|}
\hhline{-|-|-|~|-|-|-|}
{\bf Properties} &  \multicolumn{2}{l|}{$\bullet$ The {\it prefixed mesh} prior to filtering:} &  &  \multicolumn{3}{l|}{$\bullet$ Discretization error can be regulated:}  & might be less \\
  &  \multicolumn{2}{l|}{$\to$ neither adaptive nor flexible;} &  & \multicolumn{3}{l|}{$\to$ self-adaptive and variable stepsize;}  & accurate than \\
  &  \multicolumn{2}{l|}{$\to$ might fail on irregular sampling case;} &  & \multicolumn{3}{l|}{$\to$ treat irregular intervals, automatically;}   & those with \\
  &  \multicolumn{2}{l|}{$\bullet$ no control of the discretization error.} &  & \multicolumn{3}{l|}{$\bullet$ flexible to use any MATLAB ODE solver.}  & GE control.\\
\hline
\end{tabular}
}
\end{table*}

We can distinguish two basic approaches routinely being adopted in research devoted to the development of the continuous-discrete filtering methods~\cite{2012:Frogerais,2014:Kulikov:IEEE}. The first one implies the use of numerical schemes for solving the given {\it stochastic differential equation} (SDE) of the system at hand. It is usually done by using either the Euler-Maruyama method or the higher order methods based on the It\^{o}-Taylor expansion~\cite{1999:Kloeden:book}. The left panel of Table~\ref{Tab:survay} collects the continuous-discrete UKF algorithms of that type, which are already developed in engineering literature. We stress that the key property of the implementation methodology based on the SDE solvers is that the users are requested to preassign the $L$-step {\it equidistant mesh} to be applied on each sampling interval $[t_k, t_{k+1}]$ to perform the discretization, i.e. the step size is fixed. This {\it prefixed} integer $L > 0$ should be chosen prior to filtering and with no information about the discretization error arisen. Certainly, this strategy yields the UKF implementation methods with a known computational load but with no information about accuracy of the filtering algorithms. In other words, the implementation framework discussed does not ensure a good estimation quality and might fail due to a high discretization error occurred. For instance, if some measurements are missing while the filtering process, then the pre-defined $L$-step equidistant mesh might be not enough for accurate integration on the longer sampling intervals related to the missing data and this destroys the filtering algorithm. Finally, the discretization error arisen in any SDE solver is a random variable and, hence, it is uncontrollable.

An alternative implementation framework assumes the derivation of the related filters' moment differential equations (MDEs) and then an utilization of the numerical methods derived for solving ordinary differential equations (ODEs); see the right panel of Table~\ref{Tab:survay}. In~\cite{2007:Sarkka}, the MDEs are derived for continuous-discrete UKF as well as the {\it Sigma Point Differential Equations} (SPDEs), which can be solved instead. Both systems represent the ODEs that should be solve on each sampling intervals $[t_{k}, t_{k+1}]$. Clearly, one may follow the integration approach with the fixed {\it equidistant mesh} discussed above. For instance, a few steps of the Runge-Kutta methods are often suggested as the appropriate choice in engineering literature. However, this idea does not provide any additional benefit compared to the first implementation framework summarized in the left panel of Table~\ref{Tab:survay}. A rational point is to take the advantage of using ODE solvers. In particular, the discretization error arisen is now possible to control and to bound, which yields the accurate implementation way. Furthermore, the modern ODE solvers are self-adaptive algorithms, i.e. they generate the {\it adaptive} integration mesh in automatic way depending on the discretization error control involved for keeping the error within the prefixed tolerance. Such smart implementation strategy provides the accurate, self-adaptive and flexible UKF implementation methods where the users are requested to fix the tolerance value, only. In our previous works, we have proposed the continuous-discrete UKF methods by utilizing the variable-stepsize nested implicit Runge-Kutta (NIRK) formulas with the global error (GE) control suggested in~\cite{2013:Kulikov:IMA}. The main drawback is that the GE control typically yields  computationally costly algorithms. The GE is the true numerical error between the exact and numerical solutions, whereas local error (LE) is the error committed for one step of the numerical method, only. The MATLAB's built-in ODE solvers include the LE control that makes the novel continuous-discrete UKF algorithms faster than the previously derived NIRK-based methods.
In summary, the Type~II implementation framework  does not prefix the computation load because the involved ODE solvers generate the adaptive mesh depending on the problem and the discretization error control utilized, but they solve the problems with the given accuracy requirements.

The second part of each panel in Table~\ref{Tab:survay} is focused on the {\it square-root} strategies established in engineering literature for deriving the numerically stable (with respect to roundoff errors) filtering algorithms with the KF-like structure. The square-root (SR) methods ensures the theoretical properties of the filter covariance matrix in a finite precision arithmetics, which are the symmetric form and positive definiteness~\cite[Chapter~7]{2015:Grewal:book}. A variety of the SR filtering methods comes from the chosen factorization $P_{k|k}=S_{k|k}S_{k|k}^{\top}$. The traditional SR algorithms are derived by using the Cholesky decomposition, meanwhile the most recent advances are based on the singular value decomposition (SVD) as shown in~\cite{2020:Automatica:Kulikova,1986:Oshman}. For the sampled-data estimators (e.g. the CKF, GHQF and UKF), the derivation of the related SR implementations is of special interest because they demand the factorization $P_{k|k}=S_{k|k}S_{k|k}^{\top}$ in each iterate for generating the sigma/cubature/quadrature vectors.

The problem of deriving the SR methods for the UKF estimator lies in the possibly negative sigma weights used for the mean and covariance approximation. This prevents the square-root operation and, as a result, the derivation of the SR UKF algorithms. To overcome this obstacle, the {\it one-rank Cholesky update} procedure has been suggested for computing the SR of covariance matrices every time when the negative weight coefficients appear as explained in~\cite{2001:Merwe}. However, if the downdated matrix is not positive definite, then the one-rank Cholesky update is unfeasible and its failure yields the UKF estimator shutoff again. In other words, such methods do not have the numerical robustness benefit of the {\it true} SR algorithms where the Cholesky procedure is performed only once. For that reason, the SR UKF methods based on the  one-rank Cholesky updates have been called the {\it pseudo-square-root} ones in~\cite[p.~1262]{2009:Haykin}.

In this paper, we suggest the MATLAB-based SR and pseudo-SR UKF methods within both the MDE- and SPDE-based implementation frameworks as indicated in the right panel of Table~\ref{Tab:survay}.  In summary, the novel UKF methods have the following properties: (i) they are accurate due to the automatic discretization error control involved for solving the filter's expectation and covariance differential equations, (ii) they are easy to implement because of the built-in fashion of the MATLAB ODE solver in use, (iii) they are self-adaptive because the discretization mesh at every prediction step is generated automatically by the chosen ODE solver according to the discretization error control rule implemented by the MATLAB ODE solver, (iv) they are flexible, that is, any other MATLAB ODE solver of interest can be easily implemented together with its automatic discretization error control in use, and (v) they are numerically stable with respect to roundoff due to the square-root implementation fashions derived in this paper. Thus, we provide practitioners with a diversity of algorithms giving a fair possibility for choosing any of them depending on a real-world application and requirements.

Finally, it is worth noting here that the SVD-based filtering is still an open area for a future research; see Table~\ref{Tab:survay}. More precisely, the SVD solution for the Euler-Maruyama discretization-based UKF is not difficult to derive by taking into account the most recent result on the  It\^{o}-Taylor expansion-based UKF in~\cite{KuKu21IEEE_TAC}. In contrast, the derivation of the SVD-based algorithms for the Type~II framework is complicated. It demands both the MDEs and SPDEs to be re-derived in terms of the SVD SR factors instead of the Cholesky ones proposed in~\cite{2007:Sarkka}.
This is an open problem for a future research.

\section{Continuous-discrete Unscented Kalman filter}\label{problem:statement}

Consider continuous-discrete stochastic system
\begin{align}
dx(t) & = f\bigl(t,x(t)\bigr)dt+Gd\beta(t), \quad t>0,  \label{eq1.1} \\
z_k   & =  h(k,x(t_{k}))+v_k, \quad k =1,2,\ldots \label{eq1.2}
\end{align}
where  $x(t)$ is the $n$-dimensional unknown state vector to be estimated and  $f:\mathbb R\times\mathbb
R^{n}\to\mathbb R^{n} $ is the time-variant drift function. The process uncertainty is modelled by the additive noise term where
$G \in \mathbb R^{n\times q}$ is the time-invariant diffusion matrix and $\beta(t)$ is the $q$-dimensional Brownian motion whose increment $d\beta(t)$ is Gaussian white process independent of $x(t)$ and has the covariance $Q\,dt>0$. Finally, the $m$-dimensional measurement vector $z_k = z(t_{k})$ comes at some discrete-time points $t_k$ with the sampling rate (sampling period) $\Delta_k=t_{k}-t_{k-1}$. The measurement noise term $v_k$ is assumed to be a white Gaussian noise with the zero mean and known covariance $R_k>0$, $R \in \mathbb R^{m\times m}$. Finally, the initial state $x(t_0)$ and the noise processes are assumed to be statistically independent, and $x(t_0) \sim {\mathcal N}(\bar x_0,\Pi_0)$, $\Pi_0 > 0$.

The key idea behind the UKF estimation approach is the concept of Unscented Transform (UT). Following~\cite{2000:Julier}, the UT implies a set of $2n+1$ deterministically selected vectors called the {\it sigma points} generated by
\begin{align}
{\mathcal X}_{i} & = \hat x + S_x  \xi_i, \quad  i  = 0, \ldots, 2n \label{eq:ukf:vec} \\
\xi_0 & = \mathbf{0}, \xi_j  =  \sqrt{n+\lambda} \; e_j,   \xi_{n+j} = -\sqrt{n+\lambda} \; e_j \label{eq:ukf:points}
\end{align}
where $e_j$ ($j=1,\ldots, n$) stands for the $j$-th unit coordinate vector in ${\mathbb R}^n$. The matrix $S_x$ is a square-root factor of $P_x$, i.e. $P_x = S_xS_x^{\top}$  and it is traditionally defined by using the  Cholesky decomposition. Throughout the paper, we consider the lower triangular Cholesky factors.

Following the common UKF presentation proposed in~\cite{2001:Merwe}, three pre-defined scalars $\alpha$, $\beta$ and $\kappa$ should be given to calculate the weight coefficients as follows:
\begin{align}
w^{(m)}_0 & =\frac{\lambda}{n+\lambda}, & w^{(m)}_i & = \frac{1}{2n+2\lambda}, \label{eq:ukf:wmean} \\
w^{(c)}_0 & =\frac{\lambda}{n+\lambda}+1-\alpha^2+\beta, & w^{(c)}_i & = \frac{1}{2n+2\lambda} \label{eq:ukf:wcov}
\end{align}
where $i=1, \ldots, 2n$, $\lambda=\alpha^2(\kappa+n)-n$ and the secondary scaling parameters $\beta$ and $\kappa$ can be used for a further filter's tuning in order to match higher moments.

The time update step of the {\it continuous-discrete} UKF methods within Type~II framework presented in Table~\ref{Tab:survay} implies the numerical integration schemes for solving the MDEs on each sampling interval $[t_{k-1}, t_{k}]$ as follows~\cite{2007:Sarkka}:
{\small
\begin{align}
    \!\!\!\frac{d\hat x(t)}{dt} & = f\bigl(t,{\mathbb X}(t)\bigr) w^{(m)}, \label{UKF:MDE1} \\
    \!\!\!\frac{dP(t)}{dt}& = {\mathbb X}(t){\mathbb W}f^{\top}\bigl(t,{\mathbb X}(t)\bigr)\!  +\! f\bigl(t,{\mathbb X}(t)\bigr){\mathbb W}{\mathbb X}^{\top}(t)\! +\! GQG^{\top}\!\!\! \label{UKF:MDE2}
\end{align}}
 where ${\mathbb X}(t)$ stands for the matrix collected from the sigma points ${\mathcal X}_{i}(t)$ defined by~\eqref{eq:ukf:vec} around the mean $\hat x(t)$ through the matrix square-root $P^{1/2}(t)$. They are located by columns in the discussed matrix, i.e. ${\mathbb X}(t) = \Bigl[{\mathcal X}_{0}(t), \ldots, {\mathcal X}_{2n}(t) \Bigr]$ is of size $n\times(2n+1)$. Additionally,  the weight matrix ${\mathbb W}$ is defined as follows:
{\small
\begin{align}
w^{(m)}  & =\bigl[w^{(m)}_0,\ldots,w^{(m)}_{2n}\bigr]^\top, \; w^{(c)} =\bigl[w^{(c)}_0,\ldots,w^{(c)}_{2n}\bigr]^\top, \label{eq:w_c}\\
{\mathbb W} & = \bigl[I_{2n+1}-\mathbf{1}_{2n+1}^\top\otimes w^{(m)}\bigr]\mbox{\rm diag}\bigl\{w^{(c)}_0,\ldots, w^{(c)}_{2n}\bigr\} \nonumber \\
& \times \bigl[I_{2n+1}-\mathbf{1}_{2n+1}^\top\otimes w^{(m)}\bigr]^\top \label{eq:W_matrix}
\end{align}}
where the vector ${\mathbf 1}_{2n+1}$ is the unitary column of size $2n+1$ and the symbol $\otimes$ is the Kronecker tensor product.

Alternatively, one may solve the system of the UKF SPDEs derived in~\cite{2007:Sarkka}:
\begin{align}
& \frac{d{\mathbb X}^{\prime}_i(t)}{dt}  = f\bigl(t,{\mathbb X}(t)\bigr) w^{(m)} + \sqrt{n+\lambda}  \nonumber \\
&\times \bigl[{\bf 0}, \: P^{1/2}(t)\Phi(M(t)), \: -P^{1/2}(t)\Phi(M(t))\bigr]_i \label{eq2.9}
\end{align}
where the subscript $i$ refers to the $i$-th column in the related matrices, $i=0, \ldots 2n$. The notation ${\bf 0}$ stands for the zero column. The matrix $M(t)$ is computed by
\begin{align}
M(t) & =P^{-1/2}(t)\left[{\mathbb X}(t){\mathbb W} f^{\top}\bigl(t,{\mathbb X}(t)\bigr) \right. \nonumber \\
 & \left. + f\bigl(t,{\mathbb X}(t)\bigr) {\mathbb W}{\mathbb X}^{\top}(t) + GQG^\top\right]P^{-\top/2}(t) \label{eq2.10}
\end{align}
with the mapping $\Phi(\cdot)$ that returns a lower triangular matrix defined as follows: (1) split the argument matrix $M$ as  $M=\bar L + D + \bar U$ where $\bar L$ and $\bar U$ are, respectively, a strictly lower and upper triangular parts of $M$, and $D$ is its main diagonal; (2) compute $\Phi(M) = \bar L + 0.5 D$ for any argument matrix $M$.

Next, the UKF measurement update step at any time instance $t_k$ implies the generation of sigma points ${\mathcal X}_{i,k|k-1}$, $i=0,\ldots, 2n$, around the predicted mean $\hat x_{k|k-1}$ by formula~\eqref{eq:ukf:vec}  with the square-root $P_{k|k-1}^{1/2}$ of the predicted covariance matrix $P_{k|k-1}$ together with the weights in~\eqref{eq:ukf:wmean}, \eqref{eq:ukf:wcov}. For system~\eqref{eq1.1}, \eqref{eq1.2}, this allows for the following effective calculation~\cite[p.~1638]{2007:Sarkka}:
\begin{align}
{\mathbb  Z}_{k|k-1} & =h\bigl(k,{\mathbb  X}_{k|k-1}\bigr), \quad \hat z_{k|k-1} ={\mathbb  Z}_{k|k-1}w^{(m)}, \label{ckf:zpred} \\
R_{e,k} & ={\mathbb Z}_{k|k-1}{\mathbb W}{\mathbb Z}_{k|k-1}^{\top}+R_k, \label{ckf:rek} \\
 P_{xz,k} & ={\mathbb X}_{k|k-1}{\mathbb W}{\mathbb Z}_{k|k-1}^{\top}. \label{ckf:pxy} \\
 \hat x_{k|k} & =\hat x_{k|k-1}+{K}_k(z_k-\hat z_{k|k-1}),  \label{ckf:state}\\
{K}_{k} & =P_{xz,k}R_{e,k}^{-1}, P_{k|k}  = P_{k|k-1} - {K}_k R_{e,k} {K}_k^{\top} \label{ckf:gain}
\end{align}

\begin{table*}[ht!]
{\scriptsize
\renewcommand{\arraystretch}{1.3}
\caption{The {\it conventional} continuous-discrete MATLAB-based UKF filtering methods.} \label{Tab:1}
\centering
\begin{tabular}{l||l|l}
\hline
& \cellcolor{myGray} {\bf MDE-based UKF: Algorithm~1} & \cellcolor{myGray} {\bf SPDE-based UKF: Algorithm~2} \\
\hline
\hline
\textsc{Initialization:}  &  \multicolumn{2}{l}{0. Set initials $\hat x_{0|0} = \bar x_0$, $P_{0|0} = \Pi_0$. Define $\xi_i$, $i=0,\ldots, 2n$ by~\eqref{eq:ukf:points} and weights in~\eqref{eq:w_c}, \eqref{eq:W_matrix}. Set the options in~\eqref{eq2.136}.  }\\
\hline
\textsc{Time} & \multicolumn{2}{l}{1. Cholesky dec.: $P_{k-1|k-1}=P_{k-1|k-1}^{1/2}P_{k-1|k-1}^{\top/2}$. Generate ${\mathcal X}_{i,k-1|k-1}=\hat x_{k-1|k-1}+P_{k-1|k-1}^{1/2}\xi_i$, $i=0, \ldots, 2n$.} \\
\textsc{Update (TU):} & 2. Set $XP_{k-1|k-1} = [\hat x_{k-1|k-1}, P_{k-1|k-1}]$. & 2. Set ${\mathbb  X}_{k-1|k-1}=\bigl[{\mathcal X}_{0,k-1|k-1},\ldots,{\mathcal X}_{2n,k-1|k-1}\bigr]$.\\
& 3. Reshape $x^{(0)}_{k-1} = XP_{k-1|k-1}\verb"(:)"$. & 5. $\widetilde{{\mathbb  X}}^{(0)}_{k-1|k-1} = {\mathbb  X}_{k-1|k-1}\verb"(:)"$. \\
& 4. Integrate $x_{k|k-1}\leftarrow \texttt{odesolver[MDEs},x^{(0)}_{k-1},[t_{k-1},t_k]]$. & 4. $\widetilde{{\mathbb  X}}_{k|k-1}\leftarrow \texttt{odesolver[SPDEs},\widetilde{{\mathbb  X}}^{(0)}_{k-1|k-1},[t_{k-1},t_k]]$. \\
& 5. $XP_{k|k-1} \leftarrow \texttt{reshape}(x_{k|k-1}^{\texttt{end}},n,n+1)$.  & 5. ${\mathbb  X}_{k|k-1} \leftarrow \texttt{reshape}(\widetilde{{\mathbb  X}}_{k|k-1}^{\texttt{end}},n,2n+1)$. \\
& 6. Recover $\hat x_{k|k-1} = XP_{k|k-1}$\texttt{(:,1)}. & 6. Recover $\hat x_{k|k-1} = {\mathcal X}_{0,k|k-1}=[{\mathbb  X}_{k|k-1}]_1$. \\
& 7. Recover $P_{k|k-1} = XP_{k|k-1}$\texttt{(:,2:n+1)}. & 7. $P^{1/2}_{k|k-1} = \texttt{tril}([{\mathbb  X}_{k|k-1}]_{2:n+1}-\hat x_{k|k-1})/\sqrt{n+\lambda}$. \\
\hline
\textsc{Measurement} & 7a. Cholesky decomposition $P_{k|k-1}=P_{k|k-1}^{1/2}P_{k|k-1}^{\top/2}$ & $-$ Square-root $P^{1/2}_{k|k-1}$ is already available from TU. \\
\textsc{Update (MU):}  & 7b. Get ${\mathcal X}_{i,k|k-1}=\hat x_{k|k-1}+P_{k|k-1}^{1/2}\xi_i$, $i=0,\ldots,2n$. & $-$ Sigma vectors ${\mathcal X}_{i,k|k-1}$ are already defined from TU. \\
& 7c. Form ${\mathbb  X}_{k|k-1}=\bigl[{\mathcal X}_{0,k|k-1},\ldots,{\mathcal X}_{2n,k|k-1}\bigr]$. & $-$ The matrix ${\mathbb  X}_{k|k-1}$ is already available from TU. \\
\cline{2-3}
& \multicolumn{2}{l}{8. Propagate ${\mathcal Z}_{i,k|k-1}=h\bigl(k,{\mathcal  X}_{i,k|k-1}\bigr)$, $i=0,\ldots,2n$. Set ${\mathbb  Z}_{k|k-1}=\bigl[{\mathcal Z}_{0,k|k-1},\ldots,{\mathcal Z}_{2n,k|k-1} \bigr]$ of size $m\times (2n+1)$. } \\
& \multicolumn{2}{l}{9. Find $\hat z_{k|k-1}={\mathbb  Z}_{k|k-1}w^{(m)}$, $R_{e,k}={\mathbb Z}_{k|k-1}{\mathbb W}{\mathbb Z}_{k|k-1}^{\top}+R_k$, $P_{xz,k}={\mathbb X}_{k|k-1}{\mathbb W}{\mathbb Z}_{k|k-1}^{\top}$ and  ${K}_{k}=P_{xz,k}R_{e,k}^{-1}$.} \\
& \multicolumn{2}{l}{10. Update the state $\hat x_{k|k}=\hat x_{k|k-1}+{K}_k(z_k-\hat z_{k|k-1})$ and the covariance $P_{k|k}=P_{k|k-1} - {K}_k R_{e,k} {K}_k^{\top}$.} \\
\hline
\hline
Auxiliary & $[\tilde x(t)] \leftarrow \proc{MDEs}(x(t),t,n,\lambda,{\mathbb W},w^{(m)},G,Q)$
&
$[\widetilde{{\mathbb  X}}(t)] \leftarrow \proc{SPDEs}(\widetilde{{\mathbb  X}}(t),t,n,\lambda,{\mathbb W},w^{(m)},G,Q)$\\
Functions & Get matrix $X = \verb"reshape"(x,n,n+1)$; & Get matrix ${\mathbb  X}(t) \leftarrow \verb"reshape"(\widetilde{{\mathbb  X}}(t),n,2n+1)$;\\
& Recover state $\hat x(t) = [X]_1$;            & Recover $\hat x(t) = [{\mathbb  X}(t)]_{1}$;\\
& Recover covariance $P(t) = [X]_{2~:~n+1}$;  & Recover $P^{1/2}(t) = \verb"tril"\bigl([{\mathbb  X}(t)]_{2:n+1}-\hat x(t)\bigr)/\sqrt{n+\lambda}$;\\
& Factorize $P(t)$, generate nodes~\eqref{eq:ukf:points}, \eqref{eq:ukf:vec} and get ${\mathbb X}(t)$;
& Propagate $f\bigl(t,{\mathbb X}(t)$ and find ${\mathbb X}(t){\mathbb W} f^{\top}\bigl(t,{\mathbb X}(t)\bigr)$;\\
& Propagate $f\bigl(t,{\mathbb X}(t)$ and find ${\mathbb X}(t){\mathbb W} f^{\top}\bigl(t,{\mathbb X}(t)\bigr)$; & Find $M(t)$ by~\eqref{eq2.10} and split $M=\bar L + D + \bar U$;\\
& Compute the right-hand side of~\eqref{UKF:MDE1}, \eqref{UKF:MDE2}; & Compute $\Phi(M) = \bar L + 0.5 D$;\\
& Get extended matrix $\tilde X = [d\hat x(t)/dt, dP(t)/dt]$; & Find the right-hand side of~\eqref{eq2.9}, i.e. get $A = d{\mathbb  X}(t)/dt$;\\
& Reshape into a vector form $\tilde x(t)=\tilde X\verb"(:)"$. & Reshape into a vector form $\widetilde{{\mathbb  X}}(t)=A\verb"(:)"$.\\
\hline
\end{tabular}
}
\end{table*}

To provide the general MATLAB-based UKF implementation schemes, we next denote the MATLAB's built-in ODE solver to be utilized by \verb"odesolver". This means that the users are free to choose any method from~\cite[Table~12.1]{2005:Higham:book}. The key idea of using the MATLAB's built-in ODE solvers for implementing the continuous-discrete UKF is the involved LE control that allows for bounding the discretization error. The solver creates the adaptive variable stepsize mesh in such a way that the discretization error arisen is less than the tolerance value pre-defined by user. We stress that this is done in automatic way by MATLAB's built-in functions and no extra coding is required from users except for setting the ODE solvers' options prior to filtering as follows:
\begin{equation}\label{eq2.136}\scriptsize
\texttt{options = odeset('AbsTol',LET,'RelTol',LET,'MaxStep',0.1)}
\end{equation}
where the parameters \texttt{AbsTol} and \texttt{RelTol} determine portions of the {\em absolute} and {\em relative} LE utilized in the built-in control mechanization, respectively, and $0.1$ limits the maximum step size $\tau^{\rm max}$ for numerical stability reasons. Formula~\eqref{eq2.136} implies that \texttt{AbsTol=RelTol=LET} where the parameter \texttt{LET} sets the requested local accuracy of numerical integration with the MATLAB code.

The general MATLAB-based continuous-discrete UKF strategies are proposed within both the MDE and SPDE approaches in Table~\ref{Tab:1}. Since the MATLAB's built-in ODE solvers are vector-functions, one should re-arrange both the MDEs in~\eqref{UKF:MDE1}, \eqref{UKF:MDE2} and SPDEs in~\eqref{eq2.9} in the form of unique vector of functions, which is to be sent to the ODE solver. The MATLAB built-in function \texttt{reshape} performs this operation. More precisely, \texttt{A(:)} returns a single column vector of size $M\times N$ collected from columns of the given array $A \in {\mathbb R}^{N\times M}$.  Next, the built-in function \texttt{tril(A)} extracts a lower triangular part from an array $A$. Finally, the notation $\texttt{[A]}_i$ stands for the i-th column of any matrix $A$, meanwhile $\texttt{[A]}_{i:j}$ means the matrix collected from the columns of $A$ taken from the $i$th column up to the $j$th one.

We briefly discuss Algorithms~1, 2 and remark that they are of {\it conventional}-type implementations since the entire error covariance matrix, $P_{k|k}$, is updated. As a result, the Cholesky decomposition is required at each filtering step for generating the sigma vectors. However, the SPDE-based implementation in Algorithm~2 demands one less Cholesky factorization than the MDE-based scheme in Algorithm~1. Indeed, both methods imply the Cholesky decomposition at the time update step~1 but Algorithm~2 skips the factorization at step~7 of the measurement update because the propagated sigma vectors are already available. Thus, the Cholesky decomposition is avoided at each measurement update step in Algorithm~2. In step~7 of Algorithm~2, the matrix $P^{1/2}_{k|k-1}$ is recovered by taking into account formulas~\eqref{eq:ukf:vec}, \eqref{eq:ukf:points}.

The dimension of the system to be solved in Algorithm~1 is approximately two times less than in Algorithm~2. The MDEs in~\eqref{UKF:MDE1}, \eqref{UKF:MDE2} consist of $n^2+n$ equations, meanwhile the SPDEs in~\eqref{eq2.9} contain $(2n+1)n$ equations to be solved. This means that the MDE-based UKF implementation is faster than the related SPDE-based Algorithm~2, although the latter skips one Cholesky factorization at each iterate. Finally, the auxiliary functions summarized in Table~\ref{Tab:1} are intended for computing the right-hand side functions in~\eqref{UKF:MDE1}, \eqref{UKF:MDE2} and~\eqref{eq2.9}, respectively. They are presented in a vector form required by any built-in MATLAB ODEs integration scheme. It should be stressed that the MDE-based approach for implementing a continuous-discrete UKF technique demands the Cholesky factorization in each iterate of the auxiliary function for computing the sigma nodes and then for calculating the right-hand side expressions in formulas~\eqref{UKF:MDE1}, \eqref{UKF:MDE2}. This makes the MDE-based implementation strategy vulnerable
to roundoff. The filtering process is interrupted when the Cholesky factorization is unfeasible. We conclude that the MDE-based implementation way is much faster but it is also less numerically stable compared to the SPDE-based approach. Recall, the SPDE- MATLAB-based UKF implementation method controls the discretization error arisen in $(2n+1)n$ equations. This slows down the algorithm and might yield unaffordable execution time.

\section{SR solution with one-rank Cholesky updates}\label{sect3b}

Historically, the first SR methods suggested for the UKF estimator have utilized the {\it one-rank Cholesky update} procedure for treating the negative entries appeared in the coefficient vector $w^{(c)}$ and the related formulas for computing the covariance matrices. More precisely, the UKF SR solution has been developed for a particular UKF parametrization with $\alpha=1$, $\beta=0$ and $\kappa=3-n$ in~\cite{2001:Merwe}. This set of parameters yields the negative sigma points' weights $w^{(m)}_0$ and $w^{(c)}_0$ when the number of states to be estimated is $n>3$. Following the cited works, the one-rank Cholesky update is utilized as follows. Assume that $\tilde S$ is the original lower triangular Cholesky factor of matrix $\tilde P$. We are required to find the Cholesky factor of the one-rank updated matrix $P = \tilde P \pm uu^{\top}$. Given $\tilde S$ and vector $u$, the MATLAB's built-in function $S = \texttt{cholupdate}\bigl\{\tilde S,u,\pm \bigr\}$ finds the Cholesky factor $S$ of the matrix $P$, which we are looking for. If $u$ is a matrix, then the result is $p$ consecutive updates of the Cholesky factor using the $p$ columns of $u$. An error message reports when the downdated matrix is not positive definite and the failure yields the UKF estimator shutoff. Following~\cite{2009:Haykin}, the SR UKF algorithms developed within the one-rank Cholesky update procedure  are called the {\it pseudo-SR} methods.

In this paper, we suggest a general {\it array} approach for implementing the pseudo-SR UKF methods. The array form implies that the information available is collected into the unique pre-array and, next, the QR factorization is performed to obtain the post-array. All SR factors required are then simply read-off from the obtained post-array. Such algorithms are very effective in MATLAB because of vectorized operations. The key idea of our solution is to apply the QR decomposition to the part of the pre-array that corresponds to the positive entries in coefficient vector $w^{(c)}$. Next, the one-rank Cholesky updates are utilized for updating the Cholesky factors obtained by using the rows related to the negative entries in $w^{(c)}$. We illustrate our solution by the UKF parametrization with $\alpha=1$, $\beta=0$ and $\kappa=3-n$ as suggested in~\cite{2001:Merwe}. This approach can be extended in a proper way on other UKF parametrization variants.

\begin{table*}[ht!]
{\scriptsize
\renewcommand{\arraystretch}{1.3}
\caption{The pseudo-SR MATLAB-based UKF algorithms for the classical UKF parametrization with $\alpha=1$, $\beta=0$ and $\kappa=3-n$.} \label{Tab:2}
\centering
\begin{tabular}{l||l|l}
\hline
& \cellcolor{myGray} {\bf MDE-based UKF: Algorithm~1a} & \cellcolor{myGray} {\bf SPDE-based UKF: Algorithm~2a} \\
\hline
\hline
\textsc{Initialization:}  &  \multicolumn{2}{l}{0. Cholesky decomposition: $\Pi_0 = \Pi_0^{1/2}\Pi_0^{\top/2}$, $\Pi_0^{1/2}$ is a lower triangular matrix. Define $\xi_i$, $i=0,\ldots, 2n$ by~\eqref{eq:ukf:points}.} \\
& \multicolumn{2}{l}{\phantom{0.} Set initials $\hat x_{0|0} = \bar x_0$, $P_{0|0}^{1/2} = \Pi_0^{1/2}$ and options by~\eqref{eq2.136}. Find weights in~\eqref{eq:w_c}, \eqref{eq:W_matrix} and $|{\mathbb W}|^{1/2}$, $S$ by~\eqref{class_W}, \eqref{class_S}.} \\
\cline{2-3}
\textsc{Time} & \multicolumn{2}{l}{1. Generate the sigma nodes ${\mathcal X}_{i,k-1|k-1}=\hat x_{k-1|k-1}+P_{k-1|k-1}^{1/2}\xi_i$, $i=0,\ldots,2n$.} \\
\textsc{Update (TU):} & 2. $XP_{k-1|k-1} = [\hat x_{k-1|k-1}, P^{1/2}_{k-1|k-1}]$, $x^{(0)}_{k-1} = XP_{k-1|k-1}\verb"(:)"$  & $\to$ Repeat steps~2~--~7 of Algorithm~2 \\
& 3. Integrate $x_{k|k-1}\leftarrow \texttt{odesolver[MDEs-SR},x^{(0)}_{k-1},[t_{k-1},t_k]]$. & \phantom{$\to$} summarized in Table~\ref{Tab:1}.  \\
& 4. Read-off $XP_{k|k-1} \leftarrow \texttt{reshape}(x_{k|k-1}^{\texttt{end}},n,n+1)$  & \\
& 5. Get  $\hat x_{k|k-1} = XP_{k|k-1}$\texttt{(:,1)}, $P_{k|k-1}^{1/2} = XP_{k|k-1}$\texttt{(:,2:n+1)}. &  \\
\cline{2-3}
\textsc{Measurement:}  & 6. Generate ${\mathcal X}_{i,k|k-1}=\hat x_{k|k-1}+P_{k|k-1}^{1/2}\xi_i$, $i=0,\ldots,2n$. & $-$ Sigma vectors ${\mathcal X}_{i,k|k-1}$ are defined from TU. \\
\textsc{Update (MU):}  & \multicolumn{2}{l}{7. Transform ${\mathcal Z}_{i,k|k-1}=h\bigl(k,{\mathcal X}_{i,k|k-1}\bigr)$. Set ${\mathbb  Z}_{k|k-1}=\bigl[{\mathcal Z}_{0,k|k-1},\ldots,{\mathcal Z}_{2n,k|k-1} \bigr]$. Find $\hat z_{k|k-1}={\mathbb  Z}_{k|k-1}w^{(m)}$.} \\
& \multicolumn{2}{l}{8. Build ${\mathbb A}_k =
\begin{bmatrix}
[{\mathbb  Z}_{k|k-1}|{\mathbb W}|^{1/2}]_{2:2n+1}  & R_k^{1/2}
\end{bmatrix}$. Factorize ${\mathbb R}_k \leftarrow  \mbox{\texttt{qr}}({\mathbb A}_k^{\top})$. Read-off main $m\times m$ block $\tilde S_{R_{e,k}} = [{\mathbb R}_k]_m$. } \\
& \multicolumn{2}{l}{9. Update $S_{R_{e,k}} = \mbox{\texttt{cholupdate}}(\tilde S_{R_{e,k}},[{\mathbb  Z}_{k|k-1}|{\mathbb W}|^{1/2}]_{1},S_{1,1})$ and  $R_{e,k}^{1/2} = S_{R_{e,k}}^{\top}$.}\\
& \multicolumn{2}{l}{10. Find $P_{x,z} = {\mathbb X}_{k|k-1}{\mathbb W}{\mathbb Z}_{k|k-1}^{\top}$, ${K}_{k}=P_{x,z}R_{e,k}^{-{\top}/2}R_{e,k}^{-1/2}$ and $\hat x_{k|k}=\hat x_{k|k-1}+{K}_k(z_k-\hat z_{k|k-1})$.}\\
& \multicolumn{2}{l}{11. Apply $m$ one-rank updates $S_{k|k}= \mbox{\texttt{cholupdate}}(P_{k|k-1}^{{\top}/2},{\mathbb U},'-')$ where ${\mathbb U} = {K}_kR_{e,k}^{1/2}$ and get $P^{1/2}_{k|k} = S_{k|k}^{\top}$.} \\
\hline
\hline
Auxiliary & $[\tilde x(t)] \leftarrow \proc{MDEs-SR}(x(t),t,n,\lambda,{\mathbb W},w^{(m)},G,Q)$
&
$[\widetilde{{\mathbb  X}}(t)] \leftarrow \proc{SPDEs}(\widetilde{{\mathbb  X}}(t),t,n,\lambda,{\mathbb W},w^{(m)},G,Q)$\\
Functions & $X = \verb"reshape"(x,n,n+1)$, $\hat x(t) = [X]_1$, $P^{1/2}(t) = [X]_{2:n+1}$; & \qquad $\rightarrow$ Repeat from Table~\ref{Tab:1}.\\
& Generate nodes~\eqref{eq:ukf:points}, \eqref{eq:ukf:vec}, get ${\mathbb X}(t)$ and find ${\mathbb X}(t){\mathbb W} f^{\top}\bigl(t,{\mathbb X}(t)\bigr)$; & \\
& Find $M(t)$ by~\eqref{eq2.10} and split $M=\bar L + D + \bar U$; & \\
& Compute $d\hat x(t)/dt$ and $dP^{1/2}(t)/dt$ by~\eqref{UKF:MDE1}, \eqref{eq3.2:new}; & \\
& Get $\tilde X = [d\hat x(t)/dt, dP^{1/2}(t)/dt]$, reshape $\tilde x(t)=\tilde X\verb"(:)"$. & \\
\hline
\hline
& \cellcolor{myGray} {\bf MDE-based UKF: Algorithm~1b} & \cellcolor{myGray} {\bf SPDE-based UKF: Algorithm~2b} \\
\hline
\hline
\textsc{Initialization:}  & \qquad $\rightarrow$ Repeat from Algorithm~1a. & \qquad $\rightarrow$ Repeat from Algorithm~2a.\\
 \textsc{Time Update:}  & \qquad $\rightarrow$ Repeat from Algorithm~1a. & \qquad $\rightarrow$ Repeat from Algorithm~2a.\\
\cline{2-3}
\textsc{Measurement:}  & 6. Generate ${\mathcal X}_{i,k|k-1}=\hat x_{k|k-1}+P_{k|k-1}^{1/2}\xi_i$, $i=0,\ldots,2n$. & $-$ Sigma vectors ${\mathcal X}_{i,k|k-1}$ are defined from TU. \\
\textsc{Update (MU):} & \multicolumn{2}{l}{7. Transform ${\mathcal Z}_{i,k|k-1}=h\bigl(k,{\mathcal X}_{i,k|k-1}\bigr)$. Set ${\mathbb  Z}_{k|k-1}=\bigl[{\mathcal Z}_{0,k|k-1},\ldots,{\mathcal Z}_{2n,k|k-1} \bigr]$. Find $\hat z_{k|k-1}={\mathbb  Z}_{k|k-1}w^{(m)}$.} \\
& \multicolumn{2}{l}{8. Build pre-array ${\mathbb A}_k =
\begin{bmatrix}
[{\mathbb  Z}_{k|k-1}|{\mathbb W}|^{1/2}]_{2:2n+1}  & R_k^{1/2}
\end{bmatrix}$. Factorize ${\mathbb R}_k \leftarrow  \mbox{\texttt{qr}}({\mathbb A}_k^{\top})$. Read-off $\tilde S_{R_{e,k}} = [{\mathbb R}_k]_m$. } \\
& \multicolumn{2}{l}{9. Update $S_{R_{e,k}} = \mbox{\texttt{cholupdate}}(\tilde S_{R_{e,k}},[{\mathbb  Z}_{k|k-1}|{\mathbb W}|^{1/2}]_{1},S_{1,1})$ and  $R_{e,k}^{1/2} = S_{R_{e,k}}^{\top}$.}\\
& \multicolumn{2}{l}{10. Find $P_{x,z} = {\mathbb X}_{k|k-1}{\mathbb W}{\mathbb Z}_{k|k-1}^{\top}$, ${K}_{k}=P_{x,z}R_{e,k}^{-{\top}/2}R_{e,k}^{-1/2}$.
 Compute $\hat x_{k|k}=\hat x_{k|k-1}+{K}_k(z_k-\hat z_{k|k-1})$.}\\
& \multicolumn{2}{l}{11. Build ${\mathbb A}_k =
\begin{bmatrix}
[\left({\mathbb X}_{k|k-1}-{K}_{k}{\mathbb Z}_{k|k-1}\right)|{\mathbb W}|^{1/2}]_{2:2n+1}   & K_kR_k^{1/2}
\end{bmatrix}$. Factorize ${\mathbb R}_k \leftarrow  \mbox{\texttt{qr}}({\mathbb A}_k^{\top})$. Read-off $\tilde S_{P_{k|k}} = [{\mathbb R}_k]_n$. } \\
& \multicolumn{2}{l}{12. $S_{P_{k|k}} = \mbox{\texttt{cholupdate}}(\tilde S_{P_{k|k}},[\left({\mathbb X}_{k|k-1}-{K}_{k}{\mathbb Z}_{k|k-1}\right)|{\mathbb W}|^{1/2}]_{1},S_{1,1})$ and $P_{k|k}^{1/2} = S_{P_{k|k}}^{\top}$.}\\
\hline
\hline
& \cellcolor{myGray} {\bf MDE-based UKF: Algorithm~1c} & \cellcolor{myGray} {\bf SPDE-based UKF: Algorithm~2c} \\
\hline
\hline
\textsc{Initialization:} & \qquad $\rightarrow$ Repeat from Algorithm~1a. & \qquad $\rightarrow$ Repeat from Algorithm~2a.\\
 \textsc{Time Update:}  & \qquad $\rightarrow$ Repeat from Algorithm~1a. & \qquad $\rightarrow$ Repeat from Algorithm~2a.\\
\cline{2-3}
\textsc{Measurement:}  & 6. Generate ${\mathcal X}_{i,k|k-1}=\hat x_{k|k-1}+P_{k|k-1}^{1/2}\xi_i$, $i=0,\ldots,2n$. & $-$ Sigma vectors ${\mathcal X}_{i,k|k-1}$ are  defined from TU. \\
\textsc{Update (MU):}  & \multicolumn{2}{l}{7. ${\mathcal Z}_{i,k|k-1}=h\bigl(k,{\mathcal X}_{i,k|k-1}\bigr)$. Set ${\mathbb  Z}_{k|k-1}=\bigl[{\mathcal Z}_{0,k|k-1},\ldots,{\mathcal Z}_{2n,k|k-1} \bigr]$. Find $\hat z_{k|k-1}={\mathbb  Z}_{k|k-1}w^{(m)}$.} \\
& \multicolumn{2}{l}{8. Build pre-array ${\mathbb A}_k =
\begin{bmatrix}
[{\mathbb  Z}_{k|k-1}|{\mathbb W}|^{1/2}]_{2:2n+1}  & R_k^{1/2} \\
[{\mathbb  X}_{k|k-1}|{\mathbb W}|^{1/2}]_{2:2n+1}  & \mathbf{0}
\end{bmatrix}$. Factorize ${\mathbb R}_k \leftarrow  \mbox{\texttt{qr}}({\mathbb A}_k^{\top})$. Read-off $\tilde S_{A} = [{\mathbb R}_k]_{m+n}$. } \\
& \multicolumn{2}{l}{9. Update $S_{A} = \mbox{\texttt{cholupdate}}(\tilde S_{A},\begin{bmatrix}
[{\mathbb  Z}_{k|k-1}|{\mathbb W}|^{1/2}]_{1} \\
[{\mathbb  X}_{k|k-1}|{\mathbb W}|^{1/2}]_{1}
\end{bmatrix},S_{1,1})$. Get $S_{A}^{\top} =
\begin{bmatrix}
R_{e,k}^{1/2} & {\bf 0} \\
\bar P_{xz,k} & P^{1/2}_{k|k}
\end{bmatrix}$. Read-off $R_{e,k}^{1/2}$, $P_{k|k}^{1/2}$, $\bar P_{x,z}$.}\\
& \multicolumn{2}{l}{10.  Compute $K_k = \bar P_{x,z}R_{e,k}^{-1/2}$  and find $\hat x_{k|k}=\hat x_{k|k-1}+{K}_k(z_k-\hat z_{k|k-1})$.}\\
\hline
\end{tabular}
}
\end{table*}

Let us consider the measurement update step of the UKF proposed in~\cite[Algorithm~3.1]{2001:Merwe}. Following the cited paper, the Cholesky factor $R_{e,k}^{1/2}$ of the residual covariance $R_{e,k}$  is computed as follows:
{\scriptsize
\begin{align}
\tilde R_{e,k}^{1/2} & = \!qr\!\!\left[\sqrt{w_1^{(c)}}\left([{\mathcal Z}_{i,k|k-1}]_{1:2n} - \hat z_{k|k-1}\right), R_k^{1/2} \right]\!\!\! \label{eq:ukf:1rank:rek}\\
R_{e,k}^{1/2}  & = \mbox{\texttt{cholupdate}}\bigl(\tilde R_{e,k}^{1/2},[{\mathcal Z}_{i,k|k-1}]_{0} - \hat z_{k|k-1},w_0^{(c)}\bigr) \label{eq:ukf:1rank:rek2}
\end{align}
}
where $R_k^{1/2}$ is the upper triangular Cholesky factor of the measurement covariance matrix $R_k$ and the term $[{\mathcal Z}_{i,k|k-1}]_{1:2n}$ stands for a matrix collected from the vectors $\bigl[{\mathcal Z}_{1,k|k-1},\ldots,{\mathcal Z}_{2n,k|k-1} \bigr] \in {\mathbb R}^{m\times 2n}$. It is important to note that formula~\eqref{eq:ukf:1rank:rek} involves the scalar value $w_1^{(c)}$, only. The reason is that $w_1^{(c)}=w_2^{(c)}=\ldots =w_{2n}^{(c)}$ under the UKF parametrization examined in~\cite{2001:Merwe}. It can be easily extended on a more general case with any positive coefficients $w_i^{(c)}>0$, $i=1,\ldots 2n$, by introducing the diagonal matrix $[W]_{1:2n} = {\rm diag}\bigl\{w^{(c)}_1,\ldots, w^{(c)}_{2n}\bigr\}$ and utilizing $\sqrt{[W]_{1:2n}}$ instead of $\sqrt{w_1^{(c)}}$ in equation~\eqref{eq:ukf:1rank:rek}. Finally, the obtained factor $\tilde R_{e,k}^{1/2}$ is re-calculated by utilizing the one-rank Cholesky update in formula~\eqref{eq:ukf:1rank:rek2} because of possible negative value $w_0^{(c)}<0$ and unfeasible square root operation $\sqrt{w_0^{(c)}}$.

Following~\cite[p.~1638]{2007:Sarkka} and taking into account that
${\mathbb  Z}_{k|k-1} = \bigl[{\mathcal Z}_{0,k|k-1},\ldots,{\mathcal Z}_{2n,k|k-1} \bigr] \in {\mathbb R}^{m\times (2n+1)}$, the weight matrix ${\mathbb W}$ in~\eqref{eq:W_matrix} and matrix-form equation~\eqref{ckf:rek}, we represent formulas~\eqref{eq:ukf:1rank:rek}, \eqref{eq:ukf:1rank:rek2} in array form as follows:
{\scriptsize
\begin{align}
\tilde R_{e,k}^{1/2} & = qr\left[ [{\mathbb  Z}_{k|k-1}|{\mathbb W}|^{1/2}]_{1:2n}   \quad R_k^{1/2} \right], \label{eq:ukf:1rank:rek:matrix}\\
R_{e,k}^{1/2}  & = \mbox{\texttt{cholupdate}}\bigl(\tilde R_{e,k}^{1/2}, [{\mathbb  Z}_{k|k-1}|{\mathbb W}|^{1/2}]_{0},sgn(w_0^{(c)})\bigr) \label{eq:ukf:1rank:rek2:matrix}
\end{align}
}
where $[{\mathbb  Z}_{k|k-1}|{\mathbb W}|^{1/2}]_{1:2n}$ stands for a submatrix collected from all rows and the last $2n$ columns of the matrix product ${\mathbb  Z}_{k|k-1}|{\mathbb W}|^{1/2}$. The square-root matrix $|{\mathbb W}|^{1/2}$ and the related signature matrix\footnote{For zero entries $w_{i}^{(c)}=0$, $i=0,\ldots,2n$, set $\mbox{\rm sgn}(w_{i}^{(c)})=1$.} are defined by
\begin{align}
|{\mathbb W}|^{1/2}& = \bigl[I_{2n+1}-\mathbf{1}_{2n+1}^\top\otimes w^{(m)}\bigr]  \nonumber \\
& \times \mbox{\rm diag}\left\{\sqrt{|w_0^{(c)}|},\ldots, \sqrt{|w_{2n}^{(c)}|}\right\}, \label{class_W}\\
S & = \mbox{\rm diag}\left\{\mbox{\rm sgn}(w_{0}^{(c)}),\ldots,\mbox{\rm sgn}(w_{2n}^{(c)})\right\}. \label{class_S}
\end{align}

Following~\cite[Algorithm~3.1]{2001:Merwe}, the square-root factor of the error covariance matrix $P_{k|k}$ is calculated by applying the same approach to the second formula in equation~\eqref{ckf:gain}. This yields $m$ consecutive updates of the Cholesky factor using the $m$ columns of the matrix product ${\mathbb U} = {K}_kR_{e,k}^{1/2}$ because of the substraction in equation~\eqref{ckf:gain}; see~\cite[eqs.~(28),(29)]{2001:Merwe}:
\begin{align}
{\mathbb U}& = {K}_kR_{e,k}^{1/2}, S_{k|k} = \mbox{\texttt{cholupdate}}(S_{k|k-1},{\mathbb U},'-') \label{eq:ukf:1rank:P1}
\end{align}
where $S_{k|k-1}$ and  $S_{k|k}$ are, respectively, the upper triangular Cholesky factors of $P_{k|k-1}$ and $P_{k|k}$.

Finally, it is worth noting here that the time update step of the SR SPDE-based UKF variant coincides with the time update in the conventional SPDE-based UKF (Algorithm~2) because the sigma vectors are propagated instead of the error covariance matrix. Meanwhile, the MDEs-based UKF time update should be re-derived in terms of the Cholesky factors of the error covariance matrix. This problem is solved in~\cite[eq.~(64)]{2007:Sarkka} as follows:
\begin{equation}
    \frac{dP^{1/2}(t)}{dt} = P^{1/2}(t)\Phi\Bigl(M(t)\Bigr) \label{eq3.2:new}
\end{equation}
where $M(t)$ is defined by equation~\eqref{eq2.10} and the mapping $\Phi(\cdot)$ is discussed after that formula.

The first pseudo-SR variants obtained are summarized by Algorithms~1a and~2a in Table~\ref{Tab:2}. We may expect their poor numerical stability (with respect to roundoff errors) because of $m$ consecutive one-rank Cholesky updates required at the measurement undate step as suggested in~\cite{2001:Merwe}. Recall,  the downdated matrix should be positive definite, otherwise the filtering method fails due to unfeasible operation. The numerical robustness can be improved by reducing the number of the one-rank Cholesky updates involved. It is possible to perform by deriving a symmetric equation as shown in~\cite{KuKu21aEJCON}:
\begin{align}
P_{k|k}   & = \left[{\mathbb X}_{k|k-1}-{K}_{k}{\mathbb Z}_{k|k-1}\right]{\mathbb W}\left[{\mathbb X}_{k|k-1}-{K}_{k}{\mathbb Z}_{k|k-1}\right]^{\top} \nonumber  \\
 & \phantom{=} + {K}_{k}R_k {K}_{k}^{\top}. \label{eq:proof:pkk}
\end{align}
We can factorize formula~\eqref{eq:proof:pkk} as follows:
{\small
\begin{align}
\tilde P_{k|k}^{1/2} & = \!qr\!\!\left[ [\left({\mathbb X}_{k|k-1}-{K}_{k}{\mathbb Z}_{k|k-1}\right)|{\mathbb W}|^{1/2}]_{1:2n}, K_kR_k^{1/2} \right]\!\! \label{eq:ukf:1rank:p1:joseph}\\
P_{k|k}^{1/2} & = \mbox{\texttt{cholupdate}}(\tilde P_{k|k}^{1/2}, \nonumber \\
& [\left({\mathbb X}_{k|k-1}-{K}_{k}{\mathbb Z}_{k|k-1}\right)|{\mathbb W}|^{1/2}]_{0},sgn(w_0^{(c)})). \label{eq:ukf:1rank:p2:joseph}
\end{align}}

This yields two new pseudo-SR UKF Algorithms~1b and 2b summarized in Table~\ref{Tab:2}. They reduce the number of the one-rank Cholesky updates required at each iteration step from $m+1$ to $2$, Therefore, they are expected to be more stable with respect to roundoff errors.

Finally, one more pseudo-SR version can be derived for the measurement update step of the UKF estimator. It has the most simple {\it array} representation and allows for reducing a number of the one-rank Cholesky updates from two involved in Algorithms~1b and 2b to only one. First, we note that the system of equations~\eqref{ckf:gain}, \eqref{ckf:rek}, \eqref{ckf:pxy} can be summarized in the following equality\footnote{The equality can be proved by multiplying the pre- and post-arrays involved and comparing both sides of the formulas with the original UKF equations~\eqref{ckf:gain}, \eqref{ckf:rek}, \eqref{ckf:pxy}.}:
{\small
\vspace{-1cm}
\begin{align*}
& \begin{bmatrix}
{\mathbb  Z}_{k|k-1}|{\mathbb W}|^{1/2} & R_k^{1/2} \\
{\mathbb  X}_{k|k-1}|{\mathbb W}|^{1/2} & \mathbf{0}
\end{bmatrix}{\rm diag}\{S,I_m\}
\begin{bmatrix}
{\mathbb  Z}_{k|k-1}|{\mathbb W}|^{1/2} & R_k^{1/2} \\
{\mathbb  X}_{k|k-1}|{\mathbb W}|^{1/2} & \mathbf{0}
\end{bmatrix}^{\top} \\
&=\begin{bmatrix}
R_{e,k}^{1/2} & {\bf 0}_{m\times n} & {\bf 0}_{m\times (N-n)}\\
\bar P_{xz,k} & P^{1/2}_{k|k} & {\bf 0}_{n\times (N-n)}
\end{bmatrix}
\begin{bmatrix}
R_{e,k}^{1/2} & {\bf 0}_{m\times n} & {\bf 0}_{m\times (N-n)}\\
\bar P_{xz,k} & P^{1/2}_{k|k} & {\bf 0}_{n\times (N-n)}
\end{bmatrix}^{\top}
\end{align*}}
where $N=2n+1$ is the number of sigma points, $\bar P_{xz,k} = {K}_kR_{e,k}^{1/2} = P_{x,z}R_{e,k}^{-{\top}/2}$ is the normalized gain matrix, the square-root factor $|{\mathbb W}|^{1/2}$ and the signature matrix $S$ are defined by formulas~\eqref{class_W} and~\eqref{class_S}, respectively.

The following pseudo-SR UKF strategy within the one-rank Cholesky update methodology is derived
{\small
\vspace{-1cm}
\begin{align}
\tilde R & = qr\begin{bmatrix}
[{\mathbb  Z}_{k|k-1}|{\mathbb W}|^{1/2}]_{1:2n}  & R_k^{1/2} \\
[{\mathbb  X}_{k|k-1}|{\mathbb W}|^{1/2}]_{1:2n}  & \mathbf{0}
\end{bmatrix}^{\top}, \label{eq:ukf:1rank:p1:array}\\
R & = \mbox{\texttt{cholupdate}}(\tilde R, \begin{bmatrix}
[{\mathbb  Z}_{k|k-1}|{\mathbb W}|^{1/2}]_{0} \\
[{\mathbb  X}_{k|k-1}|{\mathbb W}|^{1/2}]_{0}
\end{bmatrix},sgn(w_0^{(c)})) \label{eq:ukf:1rank:p2:array}
\end{align}}
Algorithms~1c and 2c summarized in Table~\ref{Tab:2} are designed under this approach.

\section{SR solution with hyperbolic QR factorization}\label{sect3a}

\begin{table*}[ht!]
{\scriptsize
\renewcommand{\arraystretch}{1.3}
\caption{The {\it true} SR MATLAB-based UKF algorithms for the classical UKF parametrization with $\alpha=1$, $\beta=0$ and $\kappa=3-n$.} \label{Tab:5}
\centering
\begin{tabular}{l||l|l}
\hline
& \cellcolor{myGray} {\bf MDE-based UKF: Algorithm~1a-SR} & \cellcolor{myGray} {\bf SPDE-based UKF: Algorithm~2a-SR} \\
\hline
\hline
\textsc{Initialization:}  &  \multicolumn{2}{l}{0. Re-order $w^{(c)}=[w^{(c)}_{1},\ldots, w^{(c)}_{2n}, w^{(c)}_{0}]$, $w^{(m)}=[w^{(m)}_{1},\ldots, w^{(m)}_{2n}, w^{(m)}_{0}]$. Find $ {W}$, $|{\mathbb W}|^{1/2}$, $S$ by~\eqref{eq:W_matrix}, \eqref{class_W}, \eqref{class_S}.} \\
 & \qquad $\rightarrow$ Repeat from Algorithm~1a in Table~\ref{Tab:2}. & \qquad $\rightarrow$ Repeat from Algorithm~2a in Table~\ref{Tab:2}.\\
 \textsc{Time Update:}  & \qquad $\rightarrow$ Repeat from Algorithm~1a in Table~\ref{Tab:2}. & \qquad $\rightarrow$ Repeat from Algorithm~2a in Table~\ref{Tab:2}.\\
\cline{2-3}
\textsc{Measurement:}  & 1. Get ${\mathcal X}_{i,k|k-1}=\hat x_{k|k-1}+P_{k|k-1}^{1/2}\xi_i$, $i=0,\ldots,2n$. & $-$ Sigma vectors ${\mathcal X}_{i,k|k-1}$ are defined from TU. \\
\textsc{Update (MU):} & \multicolumn{2}{l}{2. ${\mathcal Z}_{i,k|k-1}=h\bigl(k,{\mathcal X}_{i,k|k-1}\bigr)$. Set ${\mathbb  Z}_{k|k-1}=\bigl[{\mathcal Z}_{0,k|k-1},\ldots,{\mathcal Z}_{2n,k|k-1} \bigr]$. Find $\hat z_{k|k-1}={\mathbb  Z}_{k|k-1}w^{(m)}$.} \\
& \multicolumn{2}{l}{3. Build ${\mathbb A}_k =
\begin{bmatrix}
R_k^{1/2} & {\mathbb  Z}_{k|k-1}|{\mathbb W}|^{1/2}  \\
\end{bmatrix}$, $J = {\mbox diag} \{I_m, S\}$. Upper triangulate ${\mathbb R}_k \leftarrow  \mbox{\texttt{jqr}}({\mathbb A}_k^{\top},J)$ by hyperbolic QR. } \\
& \multicolumn{2}{l}{4. Read-off $R_{e,k}^{1/2} = [{\mathbb R}^{\top}_k]_{m}$. Find $P_{x,z} = {\mathbb X}_{k|k-1}{\mathbb W}{\mathbb Z}_{k|k-1}^{\top}$, ${K}_{k}=P_{x,z}R_{e,k}^{-{\top}/2}R_{e,k}^{-1/2}$, $\hat x_{k|k}=\hat x_{k|k-1}+{K}_k(z_k-\hat z_{k|k-1})$.}\\
& \multicolumn{2}{l}{5. Build ${\mathbb A}_k =
\begin{bmatrix}
P_{k|k-1}^{1/2} & K_kR_{e,k}^{1/2}
\end{bmatrix}$, $J = {\mbox diag} \{I_n, -I_m\}$. Factorize ${\mathbb R}_k \leftarrow  \mbox{\texttt{jqr}}({\mathbb A}_k^{\top},J)$. Read-off $P_{k|k}^{1/2} = [{\mathbb R}^{\top}_k]_{n}$.} \\
\hline
\hline
& \cellcolor{myGray} {\bf MDE-based UKF: Algorithm~1b-SR} & \cellcolor{myGray} {\bf SPDE-based UKF: Algorithm~2b-SR} \\
\hline
\hline
\textsc{Initialization:}  & \qquad $\rightarrow$ Repeat from Algorithm~1a-SR. & \qquad $\rightarrow$ Repeat from Algorithm~2a-SR.\\
 \textsc{Time Update:}  & \qquad $\rightarrow$ Repeat from Algorithm~1a-SR. & \qquad $\rightarrow$ Repeat from Algorithm~2a-SR.\\
\cline{2-3}
\textsc{Measurement:}  & 1. Get ${\mathcal X}_{i,k|k-1}=\hat x_{k|k-1}+P_{k|k-1}^{1/2}\xi_i$, $i=0,\ldots,2n$. & $-$ Sigma vectors ${\mathcal X}_{i,k|k-1}$ are defined from TU. \\
\textsc{Update (MU):}  & \multicolumn{2}{l}{2. ${\mathcal Z}_{i,k|k-1}=h\bigl(k,{\mathcal X}_{i,k|k-1}\bigr)$. Set ${\mathbb  Z}_{k|k-1}=\bigl[{\mathcal Z}_{0,k|k-1},\ldots,{\mathcal Z}_{2n,k|k-1} \bigr]$. Find $\hat z_{k|k-1}={\mathbb  Z}_{k|k-1}w^{(m)}$.} \\
& \multicolumn{2}{l}{3. Build ${\mathbb A}_k =
\begin{bmatrix}
R_k^{1/2} & {\mathbb  Z}_{k|k-1}| {\mathbb W}|^{1/2}  \\
\end{bmatrix}$, $J = {\mbox diag} \{I_m, S\}$. Upper triangulate ${\mathbb R}_k \leftarrow  \mbox{\texttt{jqr}}({\mathbb A}_k^{\top},J)$ by hyperbolic QR. } \\
& \multicolumn{2}{l}{4. Read-off $R_{e,k}^{1/2} = [{\mathbb R}^{\top}_k]_{m}$. Find $P_{x,z} = {\mathbb X}_{k|k-1} {\mathbb W}{\mathbb Z}_{k|k-1}^{\top}$, ${K}_{k}=P_{x,z}R_{e,k}^{-{\top}/2}R_{e,k}^{-1/2}$, $\hat x_{k|k}=\hat x_{k|k-1}+{K}_k(z_k-\hat z_{k|k-1})$.}\\
& \multicolumn{2}{l}{5. ${\mathbb A}_k =
\begin{bmatrix}
K_kR_k^{1/2} & \left({\mathbb X}_{k|k-1}-{K}_{k}{\mathbb Z}_{k|k-1}\right)|{\mathbb W}|^{1/2}
\end{bmatrix}$, $J = {\mbox diag} \{I_m, S\}$. ${\mathbb R}_k \leftarrow  \mbox{\texttt{jqr}}({\mathbb A}_k^{\top},J)$. Read-off $P_{k|k}^{1/2} = [{\mathbb R}^{\top}_k]_{n}$.} \\
\hline
\hline
& \cellcolor{myGray} {\bf MDE-based UKF: Algorithm~1c-SR} & \cellcolor{myGray} {\bf SPDE-based UKF: Algorithm~2c-SR} \\
\hline
\hline
\textsc{Initialization:}  & \qquad $\rightarrow$ Repeat from Algorithm~1a-SR. & \qquad $\rightarrow$ Repeat from Algorithm~2a-SR.\\
 \textsc{Time Update:}  & \qquad $\rightarrow$ Repeat from Algorithm~1a-SR. & \qquad $\rightarrow$ Repeat from Algorithm~2a-SR.\\
\cline{2-3}
\textsc{Measurement:}  & 1. Get ${\mathcal X}_{i,k|k-1}=\hat x_{k|k-1}+P_{k|k-1}^{1/2}\xi_i$, $i=0,\ldots,2n$. & $-$ Sigma vectors ${\mathcal X}_{i,k|k-1}$ are defined from TU. \\
& \multicolumn{2}{l}{2. ${\mathcal Z}_{i,k|k-1}=h\bigl(k,{\mathcal X}_{i,k|k-1}\bigr)$. Set ${\mathbb  Z}_{k|k-1}=\bigl[{\mathcal Z}_{0,k|k-1},\ldots,{\mathcal Z}_{2n,k|k-1} \bigr]$. Find $\hat z_{k|k-1}={\mathbb  Z}_{k|k-1}w^{(m)}$.} \\
& \multicolumn{2}{l}{3. Build ${\mathbb A}_k =
\begin{bmatrix}
R_k^{1/2} & {\mathbb  Z}_{k|k-1}|{\mathbb W}|^{1/2}  \\
\mathbf{0} & {\mathbb  X}_{k|k-1}|{\mathbb W}|^{1/2}
\end{bmatrix}$, $J = {\mbox diag} \{I_m, S\}$. Upper triangulate ${\mathbb R}_k \leftarrow  \mbox{\texttt{jqr}}({\mathbb A}_k^{\top},J)$ by hyperbolic QR. } \\
& \multicolumn{2}{l}{4. Read-off the block $R = [{\mathbb R}_k]_{m+n}$ and
get $R^{\top} =
\begin{bmatrix}
R_{e,k}^{1/2} & {\bf 0} \\
\bar P_{xz,k} & P^{1/2}_{k|k}
\end{bmatrix}$. Read-off $R_{e,k}^{1/2}$, $P_{k|k}^{1/2}$ and $\bar P_{x,z}$. }\\
& \multicolumn{2}{l}{5. Compute $K_k = \bar P_{x,z}R_{e,k}^{-1/2}$  and $\hat x_{k|k}=\hat x_{k|k-1}+{K}_k(z_k-\hat z_{k|k-1})$.}\\
\hline
\end{tabular}
}
\end{table*}

In contrast to the {\it pseudo}-SR algorithms that might be unstable because of the one-rank Cholesky update procedure involved in each iterate, we may design the {\it true} SR methods. There, the Cholesky decomposition is performed only once, i.e. for the initial matrix $\Pi_0>0$. The true SR solution implies the utilization of the so-called hyperbolic QR transformations instead of the usual QR factorization, which is used for computing the Cholesky factor of a positive definite matrix given by formula $C = A^{\top}A+B^{\top}B$. In general, the UKF formulas obey the equations of the form $C = A^{\top}A \pm B^{\top}B$. The exact form depends on the UKF parametrization implemented and, in particular, on the number of negative sigma-point coefficients $w^{(c)}$ utilized. Our illustrative UKF example used for designing the one-rank Cholesky update algorithms in Table~3 implies $\alpha=1$, $\beta=0$ and $\kappa=3-n$. This yields $w^{(c)}_0<0$ when $n>3$. Taking into account this scenario, we develop and explain a general MATLAB-based SR solution within both the MDE- and SPDE-based implementation frameworks.

Following~\cite{2003:Higham}, the $J$-orthogonal matrix $Q$ is defined as one, which satisfies $Q^{\top}JQ=QJQ^{\top}=J$ where $J=\mbox{\rm diag}(\pm 1)$ is a signature matrix. The $J$-orthogonal transformations are used for computing the Cholesky factorization of a positive definite matrix $C = A^{\top}A-B^{\top}B$, where $A \in {\mathbb R}^{p\times n}$ $(p\ge n)$ and $B \in {\mathbb R}^{q\times n}$ as explained in~\cite{2003:Higham}: if we can find a $J$-orthogonal matrix $Q$ such that
\vspace{-0.5cm}
\begin{equation}
Q
\begin{bmatrix}
A \\ B
\end{bmatrix}
=
\begin{bmatrix}
R \\ 0
\end{bmatrix} \label{hyperbolic_qr}
\end{equation}
with $J = \mbox{diag}\{I_p, -I_q\}$, then $R$ is the Cholesky factor, which we are looking for. This follows from the equality
\begin{equation}
C =
\begin{bmatrix}
A \\ B
\end{bmatrix}^{\top}
J
\begin{bmatrix}
A \\ B
\end{bmatrix}
=
\begin{bmatrix}
A \\ B
\end{bmatrix}^{\top}
Q^{\top}JQ
\begin{bmatrix}
A \\ B
\end{bmatrix}
=R^{\top}R \label{proof:proof}
\end{equation}

Thus, the hyperbolic QR factorization can be used for computing the upper triangular Cholesky factor $R_{e,k}^{1/2}$ of the residual covariance $R_{e,k}$ in~\eqref{ckf:rek} as follows:
\begin{equation}
Q
\begin{bmatrix}
R_k^{1/2} \\
| {\mathbb W}|^{{\top}/2}{\mathbb  Z}_{k|k-1}^{\top}
\end{bmatrix}
=
\begin{bmatrix}
R_{e,k}^{1/2} \\ 0
\end{bmatrix} \label{hypQR:rek}
\end{equation}
where $Q$ is any $J = \mbox{diag}\{I_m, S \}$-orthogonal matrix that upper triangulates the pre-array.

We need to stress that the form of equation~\eqref{hyperbolic_qr} requires the negative weight coefficients in $w^{(c)}$ to be located at the end of the vector, otherwise the $J$-orthogonal matrix does not have the demanded form of $J = \mbox{diag}\{I_p, -I_q\}$ and the proof in~\eqref{proof:proof} does not hold. In other words, we re-order the weight vector prior to implementation as follows: $w^{(c)}=[w^{(c)}_{1},\ldots, w^{(c)}_{2n}, w^{(c)}_{0}]$ because of possible $w^{(c)}_0<0$ case. The vector $w^{(m)}$ should be re-arranged as well, i.e. $w^{(m)}=[w^{(m)}_{1},\ldots, w^{(m)}_{2n}, w^{(m)}_{0}]$. Next, the weight matrix ${\mathbb W}$ and signature $S$ are defined as usual by formulas~\eqref{class_W} and~\eqref{class_S}. Having done that, we understand that the signature matrix has the form of $S = \mbox{diag}\{{\bf 1}_{2n},-1\}$ because of the $w^{(c)}_0<0$ scenario. Hence, the $J$-orthogonal matrix in~\eqref{hypQR:rek} is $J = \mbox{diag}\{I_m, S \} = \mbox{diag}\{I_m, I_{2n},-1 \}$. Thus, we can easily prove
\begin{align*}
R_{e,k} & =
\begin{bmatrix}
R_k^{1/2} \\
| {\mathbb W}|^{{\top}/2}{\mathbb  Z}_{k|k-1}^{\top}
\end{bmatrix}^{\top}
J
\begin{bmatrix}
R_k^{1/2} \\
|{\mathbb W}|^{{\top}/2}{\mathbb  Z}_{k|k-1}^{\top}
\end{bmatrix} \\
& = R_k^{{\top}/2}I_mR_k^{1/2} + {\mathbb Z}_{k|k-1}| {\mathbb W}|^{1/2}S|{\mathbb W}|^{{\top}/2}{\mathbb Z}_{k|k-1}^{\top} \\
& =
\begin{bmatrix}
R_k^{1/2} \\
|{\mathbb W}|^{{\top}/2}{\mathbb  Z}_{k|k-1}^{\top}
\end{bmatrix}^{\top}
Q^{\top}JQ
\begin{bmatrix}
R_k^{1/2} \\
|{\mathbb W}|^{{\top}/2}{\mathbb  Z}_{k|k-1}^{\top}
\end{bmatrix} \\
& = \begin{bmatrix}
R_{e,k}^{1/2} \\ 0
\end{bmatrix}^{\top}
J
\begin{bmatrix}
R_{e,k}^{1/2} \\ 0
\end{bmatrix}
=R_{e,k}^{\top/2}I_mR_{e,k}^{1/2} + {\bf 0}S {\bf 0}^{\top}.
\end{align*}

In the same manner, an upper triangular factor of $P_{k|k}$ can be found by factorizing the second equation in~\eqref{ckf:gain}, i.e. by utilizing the hyperbolic QR factorization with  $J = \mbox{diag}\{I_n, -I_{m} \}$-orthogonal matrix as follows:
\vspace{-0.3cm} \[
Q
\begin{bmatrix}
P_{k|k-1}^{1/2} \\ R_{e,k}^{1/2} {K}_k^{\top}
\end{bmatrix}
=
\begin{bmatrix}
P_{k|k}^{1/2} \\ 0
\end{bmatrix}.
\]

\vspace{-0.3cm}
This yields the first exact SR UKF implementations with the hyperbolic QR factorizations summarized by Algorithms~1a-SR and~2a-SR in Table~\ref{Tab:5}. They are the counterparts of the one-rank Cholesky-based methods in Algorithms~1a and~2a in Table~\ref{Tab:2}. Similarly, we design the exact SR variants of Algorithms~1b and~2b and summarize them by Algorithms~1b-SR and~2b-SR in Table~\ref{Tab:5}. Finally, Algorithms~1c-SR and~2c-SR are the SR versions of the array implementations in Algorithms~1c and~2c.

\section{Numerical experiments} \label{numerical:experiments}

To investigate a difference in the numerical robustness (with respect to roundoff errors) of the suggested SR UKF methods, we explore the target tracking problem with artificial ill-conditioned measurement scheme. This scenario yields a divergence due to the singularity arisen in the covariance $R_{e,k}$ caused by roundoff~\cite[p.~288]{2015:Grewal:book}.

\begin{exampe} \label{ex:1} When performing a coordinated turn in the horizontal plane, the aircraft's dynamics obeys equation~\eqref{eq1.1} with the following drift function: $
f(\cdot)  =\left[\dot{\epsilon}, -\omega \dot{\eta}, \dot{\eta}, \omega \dot{\epsilon}, \dot{\zeta},  0, 0\right]
$ where $G={\rm diag}\left[0,\sigma_1,0,\sigma_1,0,\sigma_1,\sigma_2\right]$ with  $\omega=3^\circ/\mbox{\rm s}$, $\sigma_1=\sqrt{0.2}\mbox{ \rm m/s}$, $\sigma_2=0.007^\circ/\mbox{\rm s}$ and $Q=I_7$. The state vector is $x(t)= [\epsilon, \dot{\epsilon}, \eta, \dot{\eta}, \zeta, \dot{\zeta}, \omega]^{\top}$, where $\epsilon$, $\eta$, $\zeta$ and $\dot{\epsilon}$, $\dot{\eta}$, $\dot{\zeta}$ stand for positions and corresponding velocities in the Cartesian coordinates at time $t$, and $\omega(t)$ is the (nearly) constant turn rate. The initial $\bar x_0=[1000\,\mbox{\rm m}, 0\,\mbox{\rm m/s}, 2650\,\mbox{\rm m},150\,\mbox{\rm m/s}, 200\,\mbox{\rm m}, 0\,\mbox{\rm m/s},\omega^\circ/\mbox{\rm s}]^{\top}$ and $\Pi_0=\mbox{\rm diag}(0.01\,I_7)$.  The dynamic state is observed through the following ill-conditioned scheme~\cite{2020:Automatica:Kulikova}:
\begin{align}
\label{problem:2}
z_k & =
\begin{bmatrix}
1 & 1 & 1 & 1 & 1 &  1 &  1\\
1 & 1 & 1 & 1 & 1 &  1 &  1 +\delta
\end{bmatrix}
x_k +
\begin{bmatrix}
v_k^1 \\
v_k^2
\end{bmatrix}, \; R_k=\delta^{2}I_2
\end{align}
where parameter $\delta$ is used for simulating roundoff effect. This increasingly ill-conditioned target tracking scenario assumes that $\delta\to 0$, i.e. $\delta=10^{-1},10^{-2},\ldots,10^{-12}$.
\end{exampe}


The system is simulated on the interval $[0s, 150s]$ by Euler-Maruyama method with the step size $0.0005(s)$ for producing the exact trajectory $x^{true}_k$. For a fixed ill-conditioning parameter $\delta$ the measurements $z_k$ are defined from $x^{true}_k$ through~\eqref{problem:2} for the sampling interval $\Delta=1(s)$, $\Delta = |t_{k}-t_{k-1}|$. Next, the inverse (filtering) problem is solved by various UKF implementation methods. This yields the estimated trajectory $\hat x_{k|k}$, $k=1, \ldots, K$. The experiment is repeated for $M=100$ Monte Carlo runs for each $\delta=10^{-1},10^{-2},\ldots,10^{-12}$.  The performance of various UKF methods is assessed in the sense of the {\em Accumulated Root Mean Square Errors} in position ($\mbox{\rm ARMSE}_p$) defined in~\cite{2010:Haykin}. The MATLAB-based UKF algorithms are implemented by using the built-in ODE solver \verb"ode45" with the tolerance value $\epsilon_g = 10^{-4}$.

\begin{table}
{\tiny
\renewcommand{\arraystretch}{1.3}
\caption{The $\mbox{\rm ARMSE}_p$ (m) of the standard and SR MATLAB-based UKF methods with hyperbolic QR factorizations.} \label{tab:acc1}
\centering
\begin{tabular}{r||r|r|r||r|r|r}
\hline
$\delta$& \multicolumn{3}{c||}{\bf MDE-based UKF} & \multicolumn{3}{c}{\bf SPDE-based UKF}  \\
\cline{2-7}
& Alg.1 & Alg.1b-SR & Alg.1c-SR & Alg.2 & 2b-SR & 2c-SR  \\
\hline
    $10^{-1}$ &  7.918 &  7.918 &  7.918 & 7.921  & 7.921 & 7.921 \\ 	
    $10^{-2}$ &  {\bf fails} & 6.055 & 6.055 & 9.415  & 6.041 & 6.041 \\ 	
    $10^{-3}$ &   & 6.044 & 6.044 & {\bf fails} & 6.032 & 6.032 \\ 	
    $10^{-4}$ &   & 6.549 & 6.549 &   & 6.120 & 6.120 \\
    \hline	
    $10^{-5}$ &   & 6.068 & 6.068 &   & 7.872 & 7.872 \\ 	
    $10^{-6}$ &   & 6.069 & 6.066 &   & 7.850 & 10.19 \\ 	
    $10^{-7}$ &   & 6.066 & 6.066 &   & 7.040 & 10.19 \\ 	
    $10^{-8}$ &   & 9.216 & 6.066 &   & {\bf fails} & 10.19 \\ 	
    $10^{-9}$ &   & {\bf fails}  & 6.067 &   &  & 10.19 \\ 	
    \hline
    $10^{-10}$ &   &   &  {\bf fails} &   &  & 10.16 \\ 	
    $10^{-11}$ &   &   &  &   &  &  {\bf fails} \\ 	
\hline
\end{tabular}
}
\end{table}

The numerical stability with respect to roundoff is investigated in terms of a speed of divergence when $\delta$ tends to a machine precision limit.
Having analyzed the results presented in Table~\ref{tab:acc1}, we make a few conclusions. First, it is clearly seen that the conventional Algorithms~1 and~2 are unstable with respect to roundoff errors. They diverge fast as the problem ill-conditioning grows. The MDE-based conventional Algorithm~1 fails (due to unfeasible Cholesky factorization) when the ill-conditioning parameter $\delta=10^{-2}$. Meanwhile, the breakdown value of the SPDE-based conventional Algorithm~2 is $\delta = 10^{-3}$. This result has been anticipated and it can be explained by extra Cholesky decompositions required at each step of any MATLAB numerical scheme for computing the matrix $\mathbb X$ involved in~\eqref{UKF:MDE1}, \eqref{UKF:MDE2}. The conventional SPDE-based Algorithm~2 does not have such a drawback since it propagates the sigma vectors $\mathbb X$, straightforward. From Table~\ref{tab:acc1}, it is also clear that both the pseudo-SR and SR methods diverge much slower than the standard implementations, which process the full filter covariance matrix. Thus, they indeed improve the robustness with respect to roundoff errors, although the numerical properties of the SR methods are different. Next, we discuss them in detail.

\begin{table}
{\tiny
\renewcommand{\arraystretch}{1.3}
\caption{The $\mbox{\rm ARMSE}_p$ (m) of various pseudo-SR MATLAB-based UKF methods with one-rank Cholesky updates.} \label{tab:acc2}
\centering
\begin{tabular}{r||r|r|r||r|r|r}
\hline
$\delta$& \multicolumn{3}{c||}{\bf MDE-based UKF} & \multicolumn{3}{c}{\bf SPDE-based UKF}  \\
\cline{2-7}
& Alg.1a & Alg.1b & Alg.1c & Alg.2a & Alg.2b & Alg.2c  \\
\hline
    $10^{-1}$ &  7.918 &  7.918 &  7.918 & 7.921  & 7.921 & 7.921 \\ 	
    $10^{-2}$ &  6.056 & 6.055 & 6.055 & 9.415  & 6.041 & 6.041 \\ 	
    $10^{-3}$ &  {\bf fails} & 6.044 & 6.044 & {\bf fails} & 6.032 & 6.032 \\ 	
    $10^{-4}$ &   & 6.549 & 6.549 &   & 6.120 & 6.120 \\
    \hline	
    $10^{-5}$ &   & 6.068 & 6.068 &   & 7.874 & 7.872 \\ 	
    $10^{-6}$ &   & 6.066 & 6.066 &   & 10.16 & 10.19 \\ 	
    $10^{-7}$ &   & {\bf fails} & 6.067 &   & {\bf fails} & 10.19 \\ 	
    $10^{-8}$ &   &  & 6.066 &   &  & 10.19 \\ 	
    $10^{-9}$ &   &  & 6.067 &   &  & 10.19 \\ 	
    \hline
    $10^{-10}$ &   &   & {\bf fails} &   &  & 10.18 \\ 	
    $10^{-11}$ &   &   &  &   &  &  {\bf fails} \\ 	
\hline
\end{tabular}
}
\end{table}

 Following Table~\ref{tab:acc1}, the SPDE-based implementations are slightly more accurate than the MDE-based algorithms when $10^{-2} \ge \delta \ge 10^{-4}$, i.e. in cases of a well-conditioned scenario. Meanwhile, in the moderate ill-conditioned cases, i.e. when $10^{-5} \ge \delta \ge 10^{-8}$, the estimation errors of the SPDE-based algorithms are $1.4$ times greater than in the related MDE-based implementations. The accumulated impact of the roundoff errors is higher in case of the SPDE-based implementations compared to the MDE-based methods because of a larger ODE system to be solved. Additionally, it creates a difficulty for the MATLAB solvers for controlling the discretization error arisen and  might yield impracticable execution time.

As can be seen, the {\it array} SR Algorithms~1c-SR and~2c-SR are the most robust to roundoff. Recall, Algorithms~1b-SR and~2b-SR compute the SR factors $R_{e,k}^{1/2}$ and $P_{k|k}^{1/2}$, separately. In particular, they require two hyperbolic QR transformations at each measurement step. In contrast, Algorithms~1c-SR and~2c-SR collect one extended pre-array and apply one hyperbolic transformation for computing $R_{e,k}^{1/2}$ and $P_{k|k}^{1/2}$, in parallel. The SR UKF implementations with less hyperbolic QR transformations involved are more numerically stable.

 Let us examine the performance of the pseudo-SR UKF algorithms proposed in this paper. Following Table~\ref{tab:acc2}, a superior performance of the SPDE-based implementation way (Algorithms~2a, 2b, 2c) over the MDE-based one (Algorithms~1a, 1b, 1c) is observed in case of the well-conditioned scenario. The situation is dramatically changed in case of the moderate ill-conditioned tests where the estimation errors of the SPDE-based algorithms are greater than in the MDE-based counterparts. The impact of the accumulated roundoff errors is higher because of a larger ODE system to be solved by the SPDE-based methods. Next, the pseudo-SR approach in Algorithms~1a and~2a provides the least stable implementation way. We have anticipated this result in Section~3 because of the  $m$ consecutive one-rank Cholesky updates required. Meanwhile, Algorithms~1c and~2c are again the most stable methods, similarly to Algorithms~1c-SR and~2c-SR discussed previously. Thus, the reduced number of the one-rank Cholesky updates yields a more stable implementation strategy. The {\it array} Algorithms~1c and~2c are the most robust method among all pseudo-SR UKF variants under examination, i.e. they diverge slowly when $\delta$ tends to machine precision.

Finally, having compared the results in Tables~\ref{tab:acc1} and~\ref{tab:acc2}, we conclude that the SR approach based on the hyperbolic QR factorizations yields the same estimation quality as the pseudo-SR methodology within the one-rank Cholesky updates. This is a consequence of algebraic equivalence between the true SR methods proposed, the pseudo-SR algorithms and the standard UKF implementations (which update the full filter covariance matrix). However, it is important to stress that their stability, i.e. the speed of divergence in ill-conditioned situations is different. The most stable SR UKF implementation methods turn out to be the SPDE-based Algorithms~2c and~2c-SR (with an unique feature of being the {\it array} algorithms), although they require a two times larger ODE system to be solved than in the MDE-based UKF implementations in Algorithms~1c and~1c-SR. The users can choose any of them depending on the requirements for solving practical applications.





\end{document}